# EXISTENCE OF POLYNOMIALS WITH GIVEN ROOTS OVER NON-COMMUTATIVE RINGS

ALINA G. GOUTOR

*Abstract.* The paper studies the question of existence of polynomials with given roots over associative non-commutative rings with identity. It is shown that in the case of an associative division ring for arbitrary $n$ elements of this ring there exists a polynomial of degree $n$ whose roots are these elements. Sufficient conditions for the existence of such a polynomial are also obtained in the case of an arbitrary (not necessarily division) associative ring with identity. The case of polynomials defined over a matrix ring over a field is considered separately; for such polynomials a criterion for the existence of a second-degree polynomial with given roots is obtained; examples of constructing polynomials with given roots are also given.

*Acknowledgments.* I would like to express my gratitude to my scientific supervisor S.V. Tikhonov for interesting ideas and useful comments.

## Introduction

The paper examines the question of the existence of polynomials with given roots in polynomial rings over associative rings with identity.

Let $R$ be an associative ring with identity. We will consider polynomials of the form.

$$P(x) = a_n x^n + a_{n-1} x^{n-1} + \ldots + a_1 x + a_0, a_k \in R, \qquad (1)$$

where the variable $x$ commutes with the coefficients $a_k$. We will denote the ring of such polynomials by $R[x]$. The addition of polynomials from $R[x]$ is determined in the usual way, multiplication is defined according to the rule

$$(a_n x^n + \cdots + a_0)(b_m x^m + \cdots + b_0) = (c_{m+n} x^{m+n} + \cdots + c_0),$$

where $c_k = \sum_{i+j=k} a_i b_j$. The degree of a polynomial of the form (1) is also determined in the usual way and is equal to $n$ if $a_n \neq 0$. For $a \in R$, we define by $P(a)$ the element

$$a_n a^n + a_{n-1} a^{n-1} + \cdots + a_1 a + a_0.$$

Let us call $a \in D$ a (right) root of $P(x)$ if $P(a) = 0$. It is known that $a \in R$ is a root of a nonzero polynomial $P(x)$ if and only if $x - a$ is a right divisor of $P(x)$ in $R[x]$ [1, Prop. 16.2], i.e. $P(x) = F(x)(x - a)$ for some polynomial $F(x)$ from $R[x]$.

The question of finding the roots of polynomials over division rings is studied in ring theory and applied mathematics. The most studied is the case of polynomials with coefficients in the algebra $\mathbb{H}$ of Hamiltonian quaternions (see the works of Serodio R., Pereira E., Vitoria J., Huang L., So W., Abrate M., Chapman A., Machen C.). It is known that any polynomial from $\mathbb{H}[x]$ is decomposed into a product of linear factors. In [2], formulas for the roots of such polynomials in special cases were obtained. In [3], a formula was obtained for the roots of a quadratic polynomial in $\mathbb{H}[x]$. This formula was generalized for an arbitrary quaternion algebra in [4] and [5]. In [6] it was shown that the roots of any polynomial in $\mathbb{H}[x]$ are the roots of the so-called companion polynomial with real coefficients. In [7], an algorithm was given for finding all the roots of a polynomial in $\mathbb{H}[x]$ using companion polynomials. In [8], several of these results were generalized to the case of arbitrary division algebra. In [9], formulas were obtained for finding the roots of polynomials with coefficients in division rings in special cases.

However, the inverse problem of finding polynomials with given roots is also important. The paper [10] examines the question of the existence of polynomials with given roots in a polynomial ring over matrix rings. It was shown in [10] that for two matrices there does not always exist a quadratic polynomial that has these matrices as its roots. Note that in the case of a commutative ring, for any elements of this ring $x_1, ..., x_n$, there always exists a polynomial $(x-x_1)...(x-x_n)$ of degree $n$, the roots of which are these elements. However, in the case of non-commutative rings the situation is fundamentally different. Below we construct an example where there are no polynomials with given roots.

## The case of an arbitrary associative ring with identity

We need the following lemma about the roots of polynomials over associative rings with identity. In the case of an associative ring with division, a similar statement can be found, for example, in ([1, Pr. 16.3]). The proof of the statement for division rings can easily be generalized to the case of associative rings with identity. For the reader's convenience, we present the proof here.

*Lemma* 1. Let $R$ be an associative ring with identity and let $P(x) = L(x)Q(x) \in R[x]$. Let $d \in R$ be such that the element $h := Q(d)$ invertible. Then
$$P(d) = L(hdh^{-1})Q(d).$$

In particular, if $hdh^{-1}$ is the root of the polynomial $L(x)$, then d is the root of the polynomial $P(x)$.

*Proof.* Let $L(x) = \sum a_i x^i$, then $P(x) = \sum a_i x^i Q(x) = \sum a_i Q(x) x^i$. Then we have

$$P(d) = \sum a_i Q(d) d^i = \sum a_i h d^i h^{-1} h = \sum a_i \left(hdh^{-1}\right)^i h = L\left(hdh^{-1}\right)Q(d). \ \square$$

It follows from the previous statement that if $d$ is a root of the right factor, or if $hdh^{-1}$ (with $h$ defined above) is the root of the left factor, then $d$ is a root of the product of two polynomials. However, since there may be zero divisors in an arbitrary ring, the product may have other roots.

Let we have different $x_1, \ldots, x_n \in R$ and $R_1(x) = x - x_1$. For $i = 1, \ldots, n-1$ in the case when $R_i(x_{i+1})$ is an invertible element in $R$, recursively define

$$R_{i+1}(x) = \left(x - R_i(x_{i+1}) x_{i+1} \left(R_i(x_{i+1})\right)^{-1}\right) R_i(x).$$

Using the notation above, we obtain the following theorem.

*Theorem* 1. If for any $i = 1, \ldots, n-1$ the element $R_i(x_{i+1})$ is invertible in $R$, then $R_n(x_i) = 0$ for any $i = 1, \ldots, n$. In other words, in the case of invertibility of elements $R_i(x_{i+1})$ in the ring $R$, for any distinct elements $x_1, \ldots, x_n \in R$ there exists a polynomial in $R[x]$ of degree $n$ such that $x_1, \ldots, x_n$ are its roots.

*Proof.* Using induction, we prove that $x_1, \ldots, x_i$ are roots of $R_i(x)$, $i = 1, \ldots, n$. Indeed, when $i = 1$ the element $x_1$ will be the root of the polynomial $R_1(x) = x - x_1$. Let the statement is true for *i*. Let us prove the statement for $i + 1$, that is, we will show that $x_1, \ldots, x_{i+1}$ are the roots of the polynomial $R_{i+1}(x)$. We have

$$R_{i+1}(x) = \left(x - R_i(x_{i+1})x_{i+1}\left(R_i(x_{i+1})\right)^{-1}\right)R_i(x).$$

By assumption $x_1, \ldots, x_i$ are the roots of the right multiplier $R_i(x)$, and therefore of the entire product, i.e. of the polynomial $R_{i+1}(x)$. It remains to show that $x_{i+1}$ is also a root of the polynomial $R_{i+1}(x)$. This follows from Lemma 1. □

*Remark* 1. If under the conditions of the previous theorem $R_i(x_{i+1}) = 0$ for some *i*, then we can take $R_{i+1}(x) = R_i(x)$. If we need a polynomial of degree $i+1$, then we can take, for example, $R_{i+1}(x) = xR_i(x)$. Thus, in this case too, there is a polynomial from $R[x]$ of degree less than or equal to $n$, the roots of which will be different $x_1, \ldots, x_n \in R$.

*Theorem* 2. Let *D* be an associative division ring. Then for any elements $x_1, \ldots, x_n \in D$ there exists a polynomial $F(x)$ of degree $n$ such that $F(x_i) = 0$ for any $i = 1, \ldots, n$.

*Proof.* Follows from Theorem 1 and the previous remark. □

Thus, the construction of the polynomial from Theorem 2 occurs according to the following algorithm for finding a polynomial of degree $n$, the roots of which are different fixed elements $x_1, \ldots, x_n$. First we consider the polynomial $F_1(x) = x - x_1$, the root of which is $x_1$. Then, it is built on it $F_2(x) = (x - y_2)(x - x_1)$ such that its roots are $x_1$ and $x_2$. For this $y_2$ is taken equal to $(x_2 - x_1)x_2(x_2 - x_1)^{-1}$. For construction $F_3(x)$, the roots of which are $x_1$, $x_2$ and $x_3$, we need to check if $x_3$ is a root of $F_2(x)$. If it is true, then $F_3(x) = x\,F_2(x)$, if it is false, then $F_3(x) = (x - y_3)\,F_2(x)$, where $y_3$ is taken from Lemma 1 such that $x_3$ is a root of $F_3(x)$. Then the construction continues in a similar manner until it is built $F(x) = F_n(x)$, the roots of which will be $x_1, \ldots, x_n$.

Next, we will separately consider the case of second-degree polynomials with two given different roots $x_1$ and $x_2$.

When for specific elements $x_1$ and $x_2$ of an associative ring with identity $R$ is there a polynomial of the form

$$F(x) = x^2 + a_1x + a_0 \qquad (2)$$

the roots of which are $x_1$ and $x_2$?

Substituting first $x_1$ in (2) then $x_2$ and subtracting one equality from the other, we get

$$x_1^2 + a_1 x_1 + a_0 = 0, \ x_2^2 + a_1 x_2 + a_0 = 0, \qquad (3)$$

$$x_1^2 - x_2^2 + a_1(x_1 - x_2) = 0 \ [5, \text{p. 2}].$$

We rewrite the last equality as follows

$$a_1(x_1 - x_2) = x_2^2 - x_1^2 \qquad (4).$$

*Lemma* 1. For elements $x_1$ and $x_2$ there is a polynomial of the form (2) with roots $x_1$ and $x_2$ if and only if there is an element $a_1$ satisfying equality (4).

*Proof.* In one direction it is obvious. If there is a polynomial of the form (2) for $x_1$ and $x_2$, then $a_1$ satisfies (4). Conversely, let there exists $a_1$ satisfying equality (4) for given $x_1$ and $x_2$. This is equivalent to the fact that $x_1^2 + a_1 x_1 = x_2^2 + a_1 x_2$. If now we take $a_0$ equal to $-x_1^2 - a_1 x_1 = -x_2^2 - a_1 x_2$, then $x_1$ and $x_2$ will be the roots of the polynomial (2). □

Note that there is a polynomial of the form (2) such that $x_1$ and $x_2$ are its roots if there exists $(x_1 - x_2)^{-1}$ (this also follows from theorem 1). In this case, $a_1$ is uniquely expressed from (4), and $a_0$ is also uniquely found from one of the equalities (3). How can we check whether there is a polynomial of the form (2) whose roots are the elements $x_1$ and $x_2$ in the case where the element $x_1 - x_2$ is not invertible? When $R$ is a matrix ring over a field, this issue is considered in the next section. Note that if for some elements $x_1$ and $x_2$ there exists a quadratic polynomial of the form (2) whose roots they are, then for them there also exists a polynomial of degree $n > 2$, obtained by multiplying (2) by, for example, $x^{n-2}$.

*Theorem* 3. For elements $x_1 \neq x_2$ of an associative ring $R$ with identity a sufficient condition for the existence of an equation of degree $n$ of the form

$$x^n + a_{n-1} x^{n-1} + \ldots + a_1 x + a_0 = 0 \qquad (5)$$

of which they are solutions will be non-degeneracy of one of the elements $(x_1 - x_2)$, $\left(x_1^2 - x_2^2\right)$, …, $\left(x_1^{n-1} - x_2^{n-1}\right)$.

*Proof.* Let, for example, an element $(x_1 - x_2)$ has the inverse element. Let $a_i$, $i = 2, \ldots, n - 1$, are arbitrary elements of the ring $R$. Let us put

$$a_1 = \left(x_2^n - x_1^n\right)(x_1 - x_2)^{-1} - a_{n-1}(x_1^{n-1} - x_2^{n-1})(x_1 - x_2)^{-1} - \ldots - a_2 (x_1 - x_2)^2 (x_1 - x_2)^{-1}.$$

Then

$$a_1(x_1 - x_2) = x_2^n - x_1^n - a_{n-1}\left(x_1^{n-1} - x_2^{n-1}\right) - \ldots - a_2\left(x_1^2 - x_2^2\right).$$

Let also $a_0 = -a_1 x_1 - a_2 x_1^2 - \ldots - a_{n-1} x_1^{n-1} - x_1^n = -a_1 x_2 - a_2 x_2^2 - \ldots - a_{n-1} x_2^{n-1} - x_2^n$. Then $x_1$ and $x_2$ are the roots of the polynomial $x^n + a_{n-1} x^{n-1} + \ldots + a_1 x + a_0$. Thus, the desired equation has been found.

In case the inverse element exists for one of the elements $\left(x_1^2 - x_2^2\right)$, ..., $\left(x_1^{n-1} - x_2^{n-1}\right)$, the reasoning is similar. □

### The case of a matrix ring over a field

Let us consider the case where $R$ is a ring of square matrices over a field. As the example from [10, p. 3] shows, there exist matrices $x_1$ and $x_2$ such that there is no second-degree polynomial in $R[x]$ whose roots are these matrices $x_1$ and $x_2$. Another example of this type is given below (see Example 3). In the case of a ring of square matrices over a field, the following criteria for the existence of a second-degree polynomial with given roots are obtained.

*Theorem* 4. Let $R$ be the ring of square $k$ by $k$ matrices over a field. For $x_1$ and $x_2$ in $R$, there is a polynomial of the form (2) whose roots are them if and only if the rank of the matrix $(x_1 - x_2)^T$ equals to the rank of the matrix $\left((x_1 - x_2)^T \mid \left(x_2^T\right)^2 - \left(x_1^T\right)^2\right)$, obtained by adding rows of the matrix $\left(x_2^T\right)^2 - \left(x_1^T\right)^2$ to the right of the matrix rows $(x_1 - x_2)^T$.

*Proof.* Recall that, in view of Lemma 1, for matrices $x_1$ and $x_2$ there is a polynomial of the form (2) whose roots they are, if and only if there is a matrix $a_1$ satisfying equality (4). We transpose the equality (4):

$$(x_1 - x_2)^T a_1^T = \left(x_2^T\right)^2 - \left(x_1^T\right)^2. \tag{6}$$

From the previous equality, we note that the $j$-th element of the $i$-th column of the matrix $\left(x_2^T\right)^2 - \left(x_1^T\right)^2$ is obtained by multiplying the $j$-th row of the matrix $(x_1 - x_2)^T$ to the $i$-th column of the matrix $a_1^T$. Let us designate $a_1^T = (q_1, \ldots, q_k)$, where $q_i$ is the $i$-th column of the matrix $a_1^T$, and also similarly $\left(x_2^T\right)^2 - \left(x_1^T\right)^2 = (b_1, \ldots, b_k)$. We obtain the following equations for finding the columns of the matrix $a_1^T$:

$$(x_1 - x_2)^T q_i = b_i, \; i = 1, \ldots, k.$$

The Kronecker–Capelli theorem says that all these systems are compatible (with respect to $q_i$) if and only if the rank of the matrix $(x_1 - x_2)^T$ is equal to the ranks of the matrices $\left((x_1 - x_2)^T \mid b_i\right)$ for any $i$. Then the matrix equation (6) has a solution (there is $a_1^T$) if and only if the rank of the matrix $(x_1 - x_2)^T$ is equal to the rank of the matrix $\left((x_1 - x_2)^T \mid (x_2^T)^2 - (x_1^T)^2\right)$.

If the matrix $a_1^T$ is found, using transposition, we can also find the matrix $a_1$. Next, $a_0$ can be found from one of the equalities (3). This means that in this case there will be a polynomial of the form (2) such that $x_1$ and $x_2$ will be its roots. □

This statement can also be formulated in terms of the original matrices.

*Theorem* 5. For square $k$ by $k$ matrices $x_1$ and $x_2$, there is a polynomial of the form (2) whose roots are these matrices if and only if the rank of the matrix $x_1 - x_2$ is equal to the rank of the matrix $\begin{pmatrix} x_1 - x_2 \\ x_2^2 - x_1^2 \end{pmatrix}$, obtained by adding rows of the matrix $x_2^2 - x_1^2$ from below to the matrix rows $x_1 - x_2$.

*Proof.* In the equality (4) we denote $a_1 = \begin{pmatrix} q_1 \\ \dots \\ q_k \end{pmatrix}$, where $q_i$ are the rows of matrix $a_1$, $x_2^2 - x_1^2 = \begin{pmatrix} b_1 \\ \dots \\ b_k \end{pmatrix}$, $b_i$ are the rows of matrix $x_2^2 - x_1^2$. Then the matrix equation (4) is equivalent to the system of equations $q_i(x_1 - x_2) = b_i$, $i = 1, \dots, k$.

And the system is compatible if and only if the rank of the matrix $x_1 - x_2$ is equal to the rank of the matrix $\begin{pmatrix} x_1 - x_2 \\ b_i \end{pmatrix}$, where the last matrix is obtained by appending the row $b_i$ from below to the rows of the matrix $x_1 - x_2$. Next we obtain the required statement. □

*Example* 1. Let us give an example of such complex 2x2 matrices that there is no inverse matrix to $x_1 - x_2$, but there is a polynomial of the form (2) such that $x_1$ and $x_2$ are its roots. Such matrices will be, for example, $x_1 = \begin{pmatrix} 0 & \alpha \\ 0 & 0 \end{pmatrix}$, $x_2 = \begin{pmatrix} 0 & \beta \\ 0 & 0 \end{pmatrix}$, where $\alpha$ and $\beta$ are non-zero complex numbers, while $\alpha \neq \beta$. Since $x_1^2 = x_2^2 = \begin{pmatrix} 0 & 0 \\ 0 & 0 \end{pmatrix}$, $x_1^2 - x_2^2 = \begin{pmatrix} 0 & 0 \\ 0 & 0 \end{pmatrix}$, then the rank of the matrix $x_1 - x_2$ coincides with the rank of the matrix $\begin{pmatrix} x_1 - x_2 \\ x_2^2 - x_1^2 \end{pmatrix}$. From (4) we find that $a_1 = \begin{pmatrix} 0 & 0 \\ 0 & c \end{pmatrix}$, where $c$ is any non-zero complex number. From one of the equalities (3) we find that $a_0$ is the zero matrix. Thus, the given matrices $x_1$ and $x_2$ correspond to a continuum of polynomials of the form (2), of which they are the roots.

In the previous example, the matrices $x_1$ and $x_2$ are singular. Next, we will construct an example with non-singular matrices $x_1$ and $x_2$.

*Example* 2. Let us consider the matrices $x_1 = \begin{pmatrix} 1 & 0 \\ 0 & 1 \end{pmatrix}$, $x_2 = \begin{pmatrix} 0 & 1 \\ 1 & 0 \end{pmatrix}$. The matrix $x_1 - x_2 = \begin{pmatrix} 1 & -1 \\ -1 & 1 \end{pmatrix}$ has the rank 1. Next we get $x_1^2 = x_2^2 = \begin{pmatrix} 1 & 0 \\ 0 & 1 \end{pmatrix}$, hence $x_1^2 - x_2^2 = \begin{pmatrix} 0 & 0 \\ 0 & 0 \end{pmatrix}$.

Therefore, the rank of the matrix $x_1 - x_2$ is equal to the rank of the matrix $\begin{pmatrix} x_1 - x_2 \\ x_2^2 - x_1^2 \end{pmatrix}$ and is equal to 1. This means that there is a polynomial of the form (2) whose roots are $x_1$ and $x_2$. From (4) we find $a_1$

$$a_1 \begin{pmatrix} 1 & -1 \\ -1 & 1 \end{pmatrix} = \begin{pmatrix} 0 & 0 \\ 0 & 0 \end{pmatrix}.$$

Let $a_1 = \begin{pmatrix} a & b \\ c & d \end{pmatrix}$, then, substituting into the previous equality, we find

$$\begin{pmatrix} a & b \\ c & d \end{pmatrix} \begin{pmatrix} 1 & -1 \\ -1 & 1 \end{pmatrix} = \begin{pmatrix} a-b & -a+b \\ c-d & -c+d \end{pmatrix} = \begin{pmatrix} 0 & 0 \\ 0 & 0 \end{pmatrix}.$$

So, we get that $a_1$ is any 2 by 2 matrix that has $a = b$ and $c = d$. Let $a_1 = \begin{pmatrix} 0 & 0 \\ 0 & 0 \end{pmatrix}$, then from equalities (3) we obtain that $a_0 = \begin{pmatrix} -1 & 0 \\ 0 & -1 \end{pmatrix}$.

Note that for $x_1$ and $x_2$ there are also an infinite number of equations of the form (2) whose roots they together are.

*Proposition* 1. If for two square matrices $x_1$ and $x_2$ the equation (4) has a solution $a_1$, provided that $(x_1 - x_2)^{-1}$ does not exist, then there are infinitely many such solutions (and therefore, there are infinitely many polynomials of the form (2) with roots $x_1$ and $x_2$).

*Proof.* $a_1(x_1 - x_2) = x_2^2 - x_1^2$ Let us look at the equality (4) as a matrix equation with respect to $a_1$, then it looks like $XA = B$, where $X = a_1$, and $A$, $B$ are known matrices. Let

$$a_1 = \begin{pmatrix} q_1 \\ \ldots \\ q_k \end{pmatrix}, \text{ where } q_i \text{ are the rows of matrix } a_1, \quad x_2^2 - x_1^2 = \begin{pmatrix} b_1 \\ \ldots \\ b_k \end{pmatrix}, \quad b_i \text{ are the rows of matrix }$$

$x_2^2 - x_1^2$. The matrix equation (4) is then equivalent to the system of equations

$$q_i(x_1 - x_2) = b_i, \, i = 1, \ldots, k.$$

It is known that (4) has a solution, which means that the previous equalities also have one. It is also known that the matrix of the system $(x_1 - x_2)$ is degenerate for any $i$, which means that the previous system of equations has infinitely many solutions. And this means that (4) has infinitely many solutions as well. □

Let us give an example showing that the sufficient conditions from Theorem 3 above are not necessary for the existence of a third-degree polynomial for two matrices, and also that square matrices that do not satisfy such conditions exist.

*Example 3.* Let us consider $x_1 = \begin{pmatrix} 0 & 0 \\ 1 & -1 \end{pmatrix}$, $x_2 = \begin{pmatrix} 0 & 0 \\ 0 & 1 \end{pmatrix}$. Note that both matrices $x_1 - x_2 = \begin{pmatrix} 0 & 0 \\ 1 & -2 \end{pmatrix}$ and $x_1^2 - x_2^2 = \begin{pmatrix} 0 & 0 \\ -1 & 0 \end{pmatrix}$ are singular. This means that the sufficient conditions for the existence of a third-degree polynomial whose roots would be $x_1$ and $x_2$ are not met. There is also no equation of degree 2 for these matrices, since the rank of the matrix $\begin{pmatrix} x_1 - x_2 \\ x_2^2 - x_1^2 \end{pmatrix}$ is equal to 2, which does not coincide with the rank of the matrix $x_1 - x_2$. Let us show that, despite this, an equation of the form (5) for $n = 3$ and for $x_1$ and $x_2$ exists. Let us first substitute $x_1$ into (5) for $n = 3$, and then $x_2$, and subtract the second equality from the first

$$a_2 \begin{pmatrix} 0 & 0 \\ -1 & 0 \end{pmatrix} + a_1 \begin{pmatrix} 0 & 0 \\ 1 & -2 \end{pmatrix} = \begin{pmatrix} 0 & 0 \\ -1 & 2 \end{pmatrix}.$$

Let $a_2$ be the zero matrix, and we will search $a_1$ in the form $a_1 = \begin{pmatrix} a & b \\ c & d \end{pmatrix}$. Then

$$\begin{pmatrix} a & b \\ c & d \end{pmatrix}\begin{pmatrix} 0 & 0 \\ 1 & -2 \end{pmatrix} = \begin{pmatrix} 0 & 0 \\ -1 & 2 \end{pmatrix}.$$

From this we get that $b = 0$, $d = -1$, $a$ and $c$ are any complex numbers. Then one can take $a_1 = \begin{pmatrix} 1 & 0 \\ -1 & -1 \end{pmatrix}$. Substituting into (5) with $n = 3$ the matrices $a_i$ and $x_2$, we find the matrix $a_0$. The resulting equation of degree 3, whose roots are $x_1$ and $x_2$, has the form

$$x^3 + \begin{pmatrix} 1 & 0 \\ -1 & -1 \end{pmatrix} x = 0.$$

Let us give an example that shows that there are matrices $x_1$ and $x_2$ for which there is no polynomial of the third degree. Note that suitable matrices will be such, for example, that the first columns of the matrices $x_1 - x_2$ and $x_1^2 - x_2^2$ are zero, and in matrix $x_1^3 - x_2^3$ the first column is no longer zero.

*Example* 4. Such matrices will be $x_1 = \begin{pmatrix} 1 & -1 & 0 \\ -1 & 1 & 0 \\ 1 & 0 & 0 \end{pmatrix}$, $x_2 = \begin{pmatrix} 1 & 1 & 2 \\ -1 & 1 & 0 \\ 1 & 0 & 0 \end{pmatrix}$.

If we first substitute $x_1$ into the equality (5) for $n = 3$, then substitute $x_2$, and then subtract the second equality from the first, we get

$$a_2 \begin{pmatrix} 0 & -4 & -2 \\ 0 & 2 & 2 \\ 0 & -2 & -2 \end{pmatrix} + a_1 \begin{pmatrix} 0 & -2 & -2 \\ 0 & 0 & 0 \\ 0 & 0 & 0 \end{pmatrix} = \begin{pmatrix} -2 & 8 & 4 \\ 0 & -6 & -4 \\ 0 & 4 & 2 \end{pmatrix}.$$

The first column of the matrix on the left side of this equality consists of zeros for any matrices $a_1$ and $a_2$, but the first column of the matrix on the right side is non-zero. This means that there is no equation of the form (5) for $n = 3$, the roots of which are the matrices $x_1$ and $x_2$.

Goutor: Belarusian State University, Nezavisimosti Ave., 4, 220030, Minsk, Belarus

Email address: goutor@bsu.by